\input amstex 
\documentstyle{amsppt}
\input bull-ppt
\define \SU{\operatorname{SU}}

\keyedby{BULL301/KXJ}

\topmatter
\cvol{27}
\cvolyear{1992}
\cmonth{October}
\cyear{1992}
\cvolno{2}
\cpgs{239-242}
\ratitle
\title Smooth static solutions\\
of the Einstein-Yang/Mills equation\endtitle
\author J. Smoller, A. Wasserman, S. T. Yau, and B. 
M\lowercase{c}Leod\endauthor
\shorttitle{Einstein-Yang/Mills equations}
\address {\rm(J. Smoller and A. Wasserman)}
Department of Mathematics, 
University of Michigan, Ann Arbor, Michigan 
48109-1003\endaddress
\ml Joel\_Smoller\@ub.cc.umich.edu\newline
\indent {\it E-mail address\/}\,: 
Arthur\_Wasserman\@um.cc.umich.edu\endml
\address {\rm(B. McLeod)} Department of Mathematics,
University of Pittsburg, Pittsburgh,
Pennsylvania 15260\endaddress
\address {\rm(S. T. Yau)} Department of Mathematics,
Harvard University, Cambridge, Massachusetts
02198\endaddress
\ml styau\@math.nthu.edu.tw\endml
\date November 5, 1991 and, in revised form, January 29, 
1992\enddate
\abstract We consider the Einstein/Yang-Mills equations in 
$3+1$ space
time dimensions with $\SU(2)$ gauge group and prove 
rigorously the
existence of a globally defined smooth static solution. We 
show that
the associated Einstein metric is asymptotically flat and 
the total
mass is finite. Thus, for non-abelian gauge fields the 
Yang/Mills
repulsive force can balance the gravitational attractive 
force and
prevent the formation of singularities in 
spacetime.\endabstract
\endtopmatter

\document

\heading1\endheading

The only static, i.e., time independent, solution to the 
vacuum Einstein
equations for the gravitational field $R_{ij}-\tfrac12 
Rg_{ij}=0$ is the
celebrated Schwarzschild metric that is singular at $r=0$ 
\cite{1}.
Despite this defect, this solution has applicability for 
large $r$ to
physical problems, e.g., the perihelion shift of Mercury. 
Similarly, the
Yang/Mills equations $d^*F=0$, which unify electromagnetic 
and nuclear
forces, have no static regular solutions on $\Bbb R^4$ 
\cite{3}.
Furthermore, if one couples Einstein's equations to 
Maxwell's equations,
to unify gravity and electromagnetism\footnote""{The first 
author's
research was supported in part by NSF Contract No.\ 
DMS-89-05205 and, with the
second author, in part by ONR Contract
No.\ DOD-C-N-00014-88-K-0082; the third author was 
supported in part
by DOE Grant
No.\ DE-FG02-88ER25065; the fourth author was supported in 
part by
the U.K.\ Science and Engineering Council.}
$$R_{ij}-\tfrac12 Rg_{ij}=\sigma T_{ij},\qquad d^*F=0\tag1$$
($T_{ij}$ is the stress-energy tensor relative to the 
electromagnetic field
$F_{ij}$), the only static solution is the 
Reissner-Nordstr\"om metric,
which is again singular at the origin \cite{1}. Finally, the
Einstein-Yang/Mills (EYM) equations, which unify 
gravitational and nuclear forces,
were shown in \cite{4} to have no static regular solutions 
in $(2+1)$ space
time dimensions for any gauge group $G$. We announce here 
that the contrary
holds in $(3+1)$ space-time dimensions. Indeed, with 
$\SU(2)$ gauge group
(i.e., the weak nuclear force) we prove that the EYM 
equations (c.f.\ (1),
where now $F_{ij}$ is the $\operatorname{su}(2)$-valued 
Yang/Mills field),
admit nonsingular static solutions, whose 
metric is asymptotically flat, i.e.,
Minkowskian. (Strong numerical evidence for 
this conclusion was obtained
by Bartnik and McKinnon \cite{2} who also derived the 
relevant equations.)
Thus for non-abelian gauge fields, the Yang-Mills 
repulsive force can
balance gravitational attraction and prevent the formation 
of
singularities in spacetime. Viewed differently from a 
mathematical
perspective, it is the nonlinearity of the corresponding 
Yang/Mills
equations that allows the existence of smooth solutions.

The EYM equations are obtained by minimizing the action
$$\int(-R+|F|^2)
\sqrt g\,dx,$$
over all metrics $g_{ij}$ having signature $(-,+,+,+)$. 
These equations
become
$$R_{ij}=2F_{ik}F^k_j-\tfrac12|F|^2g_{ij}.$$
Here $R$ is the scalar curvature associated to the metric 
$g_{ij}$
and $F$ is the Yang-Mills curvature. These formidable 
equations become
more tractible if we consider static symmetric solutions.

\heading2\endheading

The problem of finding static, symmetric nonsingular 
solutions of the
EYM equations with $\SU(2)$ gauge group can be reduced to 
the study
of the following system of ordinary differential equations
$$r^2Aw''+\Phi w'+w(1-w^2)=0,\tag 2a$$
$$rA'+(2w^2+1)=1-\frac{(1-w^2)^2}{r^2},\tag 2b$$
$$2raT'+(2w\prime^2A+\Phi/ r)T=0.\tag 2c$$
Here $\Phi(r)=r(1-A)-\frac{(1-w^2)^2} r$, $A$ and $T$ are 
the unknown
metric coefficients, $ds^2=-T^{-2}(r)dt^2+A^{-1}(r)dr^2+
r^2(d\theta^2+
\sin^2\theta d\phi^2)$, and $w$ is the ``connection 
coefficient"
relative to the sought-for connection 
$\alpha=w\tau_1d\theta+[\cos
\theta\tau_3+w \sin \theta\tau_2]d\phi$, $\tau_1, \tau_2$, 
and $\tau_3$
being the generators of the Lie algebra 
$\operatorname{su}(2)$. The associated
curvature $F-d\alpha+~\alpha \wedge \alpha$ is
$$F=w'\tau_1dr\wedge d\theta+w'\tau_2dr\wedge (\sin \theta 
d\phi)
-(1-w^2)\tau_3d\theta\wedge(\sin\theta d\phi)
.$$
If $<\tau_i,\tau_j>=-2tr\tau_i\tau_j$ denotes the Killing 
form on
$\operatorname{su}(2)$, and if 
$|F|^2=g^{ij}g^{kl}F_{ij}F_{kl}$, then an easy
calculation gives
$$|F|^2=2w\prime^2/r^2+(1-w^2)^2/r^4.$$

In order that our solution has finite mass, i.e., that
$\lim_{r\to \infin}r(1-A(r))<\infin$ we require that
$$\lim_{r\to \infin} (w(r), w'(r)) \text{ be finite}.\tag3$$
Furthermore, asymptotic flatness of the metric means that
$$\lim_{r\to \infin}(A(r),T(r))=(1,1).\tag4$$
Finally, the conditions needed to ensure that our solution 
is
nonsingular at $r=0$ are
$$w(0)=1,\qquad w'(0)=0,\qquad A(0)=1,\qquad T'(0)=0.$$

One sees from (2) that the first two equations do not 
involve $T$.
Thus we first solve these for $A$ and $w$, subject to the 
above
initial and asymptotic conditions.

\heading3\endheading

We prove that under the above boundary conditions, every 
solution
is
uniquely determined by $w''(0)$; $w''(0)=-\lambda$ is a free
parameter. We seek a $\lambda>0$ such that there exists an 
orbit
$(w(r,\lambda),w'(r,\lambda))$ that ``connects two rest 
points."
It is then
not very difficult to prove that (4) will also hold.

A major difficulty is to show that the equations (2a), 
(2b) actually
define a nonsingular orbit; i.e., that $w'(r,\lambda)$ is 
bounded
and that $A(r,\lambda)$ remains positive. Our first result 
is

\thm{Theorem 1} If $0\le \lambda \le 1$, then in the region
$$\Gamma=\{w^2\le 1, w'\le 0\},$$
$A(r,\lambda)>0$ and $w'(r,\lambda)$ is bounded from below.
\ethm

On the other hand, we can also prove (see Figure 1)

\thm{Theorem 2} If $\lambda>2$, then the solution of 
equations {\rm(2a)},
{\rm(2b)},
with initial conditions {\rm(5)} blows up in $\Gamma;$ 
i.e., $w'(r)$ is
unbounded.
\ethm

\fg{11pc}\caption{\smc Figure 1}
\endfg

If $\lambda$ is near zero, then by rescaling we can show 
that the
orbit
$(w(r,\lambda),w'(r,\lambda))$ exits $\Gamma$ through the 
line $w=-1$.
Furthermore, for $\lambda=1$, numerical approximations 
indicate that
$w'$ becomes positive in the region $-1<w<0$. If this 
could be
established rigorously, we could assert the existence of 
some
$\overline{\lambda}$, $0<\overline{\lambda}<1$, for which 
the
corresponding orbit stays in $\Gamma$ for all $r\ge 0$, 
thereby
proving (3). It would then be possible to prove that
$$\lim_{r\to \infin}(w(r,\overline{\lambda})), 
w'(r,\overline{\lambda}))=
(-1,0),\tag5$$
and as a consequence, that (4) would also hold.

\heading4\endheading

We can give a completely rigorous proof of the existence 
of a connecting
orbit with $\lambda<2$, which we now outline. First 
Theorem 2 and the
fact that for $\lambda$ near 0 the corresponding orbit 
exits $\Gamma$
through the line $w=-1$ implies that there is a smallest 
$\lambda=
\overline{\lambda}$ for which the orbit 
$(w(r,\overline{\lambda})
,w'(r,\overline\lambda)
)$ does
not exit $\Gamma$ through this line. Thus only the 
following two
possibilities can arise:

$({\roman P}_1)$ There is a real number $\overline r>0$ 
such that either
(a) $w'(\overline r,\overline{\lambda})=0$, or (b) 
$A(\overline r,\overline
\lambda)=0$, or (c) $w'(r,\lambda)$ is unbounded for $r$ 
near $\overline r$.

$({\roman P}_2)$ For all $r>0$, 
$w(r,\overline{\lambda})>-1$,
$w'(r,\overline{\lambda})<0$, and 
$A(r,\overline{\lambda})>0$.

In the case that $({\roman P}_2)$ holds, we can show, as 
above, that both
(6) and (7) hold. In order to rule out possibility 
$({\roman P}_1)$, we
consider several cases. The crucial case occurs when 
$A(\overline r,
\overline{\lambda})=0$, $w'(r,\overline{\lambda})$ is 
unbounded near
$r=\overline r$, and $\Phi(\overline r, 
\overline{\lambda})=0.$ Now set
$\overline w=\lim_{r\nearrow \overline r} 
w(r,\overline{\lambda})$. If
$\overline w<0$, then defining 
$v(r,\lambda)=(Aw')(r,\lambda)$, we
show that $v$ satisfies a first order ode, and we can 
prove that for
$\lambda<\overline{\lambda}$, $\lambda$ near 
$\overline{\lambda}$,
there is an $r=r(\lambda)$ such that $v(r,\lambda)=0$ and
$w(r,\lambda)>-1$. This violates the definition of 
$\overline{\lambda}$.
Similarly, if $\overline w>0$, we can reduce this case to 
the previous
one. Finally, the case where $\overline w=0$ is dealt with 
by extending
our solution into the complex plane and using the fact 
that the pair
of functions $(w(r),A(r))=(0,1+1/r^2-c/r)$ is always a 
solution of
(2a) and (2b).

\enddocument